\documentclass[num-refs]{wiley-article}

\usepackage{siunitx}

 \newcounter{hl_crl}
 \newtheorem{crl}[hl_crl]{Corollary}

 \newcounter{hl_thm}
 \newtheorem{theo}[hl_thm]{Theorem}

\newcommand{\smf}{\setlength{\arraycolsep}{2pt}}
\newcommand{\smfs}{\small\setlength{\arraycolsep}{1pt}}

\newcommand{\vecdx}{\left[\begin{array}{l}\dot{x}_1\\\dot{x}_2\end{array}\right]}
\newcommand{\vecx}{\left[\begin{array}{l}{x}_1\\{x}_2\end{array}\right]}
\newcommand{\splitA}{\left[\begin{array}{cc}A_{11}&A_{12}\\A_{21}&A_{22}\end{array}\right]}
\newcommand{\vechx}{\left[\begin{array}{l}{x}_1\\\hat{x}_2\end{array}\right]}
\newcommand{\splitAbar}{\left[\begin{array}{cc}\bar{A}_{11}&\bar{A}_{12}\\\bar{A}_{21}&\bar{A}_{22}\end{array}\right]}
\newcommand{\splitB}{\left[\begin{array}{c}B_{1}\\0\end{array}\right]}
\newcommand{\splitD}{\left[\begin{array}{c}D_{1}\\D_{2}\end{array}\right]}
\newcommand{\norm}[1]{\| #1 \|}

\newlist{hip}{enumerate}{1}
\setlist[hip]{label=\textbf{(H.\arabic*)},labelindent=\parindent, leftmargin=*,nolistsep}
\setlist[hip]{label={(A.\arabic*)},labelindent=\parindent, leftmargin=*,nolistsep}

\papertype{Preprint submitted to arxiv.org}
\paperfield{2025}

\title{Sliding Mode Control for Uncertain Systems with Time-Varying Delays via Predictor Feedback and Super-Twisting Observer}

\author[1]{Hardy Pinto}
\author[2]{Tiago Roux Oliveira}
\author[3]{Liu Hsu}

\affil[1]{PETROBRAS---Centro de Pesquisas Leopoldo Américo Miguez de Mello, Rio de Janeiro -- RJ, Brazil.}
\affil[2]{Department of Electronics and Telecommunication Engineering (DETEL), State University of Rio de Janeiro (UERJ), Rio de Janeiro -- RJ, Brazil.}
\affil[3]{Program of Electrical Engineering (PEE), Federal University of Rio de Janeiro (UFRJ/COPPE), Rio de Janeiro -- RJ, Brazil.}

\corraddress{Corresponding author: T. R. Oliveira}
\corremail{tiagoroux@uerj.br}

\fundinginfo{This research was supported by Coordenação de Aperfeiçoamento de Pessoal de Nível Superior (CAPES), Conselho Nacional de Desenvolvimento Científico e Tecnológico (CNPq), and Fundação de Amparo à Pesquisa do Estado do Rio de Janeiro (FAPERJ).}

\runningauthor{Hardy Pinto et al.}

\begin{document}

\begin{frontmatter}
\maketitle

\begin{abstract}
This paper introduces a novel stabilization control strategy for linear time-invariant systems affected by known time-varying measurement delays and matched unknown nonlinear disturbances, which may encompass actuator faults. It is considered that part of the state vector is not available for real-time measurement. To address this, the proposed approach combines an open-loop predictor with a state observer designed using the Super-Twisting Algorithm, aiming to compensate for the delays and estimate the unmeasured state components.
Specifically, the nonlinear observer-based framework enables the reconstruction of unmodeled fault signals without assuming that they originate from a known exogenous system, offering robustness against parametric uncertainties. Meanwhile, the predictor forwards the delayed output in time. Subsequently, a sliding mode control law is formulated to enforce an ideal sliding mode and ensure global stabilization, even under a broader class of perturbations, unmodeled disturbances, parametric uncertainties, and delays, owing to the integration of the Super-Twisting observer. Numerical simulations illustrate the efficiency of the proposed approach.
\end{abstract}
\end{frontmatter}

\section{Introduction}
Delays in actuation or measurement are commonly found in various control engineering applications \cite{Krstic-2010-2}, including chemical processes, machining, combustion systems, and teleoperation. Significant measurement delays, in particular, can lead to instability and poor performance in control loops.

For systems experiencing input-output delays, control laws are typically designed using prediction feedback strategies \cite{Mirkin-2003}. In cases where the systems are stable and measurement delays are constant, the control law often incorporates the \textit{Smith Predictor} \cite{Smith-1957}. This approach has also been adapted for unstable systems or those described by state-space models \cite{Artstein-1982}. In \cite{Krstic-2010}, time-varying delays of arbitrary lengths are compensated for through a predictor-based method derived from the \textit{PDE-backstepping} technique \cite{Smyshlyaev-2005}, where delays are modeled as a partial differential equation (PDE). Similarly, \cite{Lechappe-2018} presents a Lyapunov-based predictor that can simultaneously address both input and output delays. While these methods are widely adopted, many of the proposed predictors necessitate complete knowledge of the system model and delay parameters, which limits their robustness against uncertainties or disturbances \cite{Lechappe-2015}. Adaptive techniques have been introduced to overcome some of these limitations \cite{Delphine-2009, Krstic-2010-2}. Another commonly employed approach for handling measurement delays involves the use of cascade observers \cite{Lagarrigue-paper, Camila-IJC-paper}. This method is applicable to systems with arbitrary measurement delays, although it was initially designed for constant delays \cite{Lagarrigue-paper}. However, even when delays are short in duration, the large number of observers in the cascade can make this approach difficult to implement in certain applications \cite{Holloway-2016}.

The increasing use of remotely controlled systems has spurred greater interest in fault reconstruction techniques, which enable both fault diagnosis and closed-loop control by estimating the system's state and disturbances that reflect actuator faults. Among the methods used for simultaneous state vector and unknown disturbance estimation, key approaches include unknown input observers (UIO) \cite{Patton-1993, Hui-2005, Chakrabarty-2017} and sliding mode observers (SMO) \cite{Patton-paper, Tan-2003}. In \cite{Hardy-SBAI2017,Hardy-IJRNC_bayer_mengao}, a UIO was developed to handle structured failures, utilizing the \textit{PDE-backstepping} technique \cite{Smyshlyaev-2005}, while a sliding mode control (SMC) law \cite{Utkin-book} was employed to mitigate the impact of actuator faults on the output, which was also affected by a time-varying delay. The sliding mode control (SMC) approach developed in \cite{Hardy-SBAI2017,Hardy-IJRNC_bayer_mengao} employed a disturbance observer designed for structured disturbances. This method operated under the assumption that disturbances originated from a known linear exogenous dynamic system, but did not account for potential parametric uncertainties in the system.


This work presents a new approach for the stabilization of linear time-invariant systems in which a portion of the state vector is not accessible in real-time and the output is subject to a known, time-varying measurement delay. Additionally, the proposed method enables both the rejection and the real-time estimation of nonlinear disturbances, which may represent actuator faults from a diagnostic perspective. By employing an open-loop predictor in combination with a nonlinear observer grounded on second-order sliding mode techniques—specifically the \textit{Super-Twisting Algorithm} (STA) \cite{Nagesh-2014}—it becomes feasible to reconstruct the complete state vector at the current instant and simultaneously estimate the fault affecting the actuator. A SMC strategy is then applied to ensure global stabilization of the closed-loop system, even when subjected to parametric uncertainties and unknown, possibly nonlinear, disturbances.

Moreover, the findings indicate that, provided the predictors are capable of compensating the output delays and the sliding motion is enforced in finite time regardless of the initial conditions, the proposed sliding mode controller ensures asymptotic stabilization, effective disturbance rejection, and accurate fault reconstruction. While the design of explicit predictors for uncertain nonlinear systems with delays is generally uncommon \cite{Nikos:2013}, it is feasible to integrate predictor-based feedback for uncertain linear systems subject to nonlinear disturbances with robust nonlinear control strategies, thereby broadening the scope of time-delay systems that can be effectively handled. 
A simulation example is presented to demonstrate the theoretical developments and validate the effectiveness of the proposed approach.

\section{Problem Statement}

Let us analyze a linear time-invariant system affected by actuator faults, modeled as an unknown disturbance input $d(t) \in \mathbb{R}^p$, along with parametric uncertainties captured by the term $\delta(x,t)$:
\smf\begin{align}
\left\lbrace
\begin{array}{lll}
\dot{x}(t)&=&Ax(t)+B \left[u(t)+ d(t)\right] + D \delta(x,t)\,,\\
y(t)&=&Cx(t-\tau(t)),
\end{array}
\right.
\label{eqn:system_eq}
\end{align}\normalsize
where $x(t)\in \mathbb{R}^n$ is the state vector and $u(t) \in \mathbb{R}^m$ is the control input. The output signal $y(t)\in \mathbb{R}^p$, $q \leq p <n$, has a known time-varying measurement delay, denoted by $\tau(t)$. The matrices $A \in \mathbb{R}^{n \times n}$, $B\in \mathbb{R}^{n \times m}$, $D\in \mathbb{R}^{n \times h}$ and $C\in \mathbb{R}^{p \times n}$ are also assumed to be known. 

To simplify the notation, we denote by $x := x(t)$ the state evaluated at the current time, while delayed terms are represented by $x_\tau := x(t - \tau(t))$.

In this work, we consider that the state vector $x$ is only partially available for measurement, which allows system (\ref{eqn:system_eq}) to be reformulated as:
\smfs\begin{align}
\left\lbrace
\begin{array}{lll}
\vecdx&=&\splitA\vecx+ \splitB (u+d) + \splitD \delta(x,t)\,,\\
y(t)&=& x_{\tau,2}\,,
\end{array}
\right.
\label{eqn:system_eq_split_incerto}
\end{align}\normalsize
with \( x_{\tau,2} := x_2(t - \tau(t)) \), the state component \( x_1 \in \mathbb{R}^{n-p} \) is assumed to be measurable, whereas \( x_2 \in \mathbb{R}^p \) remains unmeasured. The matrices \( A_{11},\;A_{12},\;A_{21},\;A_{22} \), and \( B_1 \) are of appropriate dimensions. In this framework, \( x_1 \) may characterize the dynamics of actuator subsystems—for example, pilot valves in oil industry applications—while \( x_2 \) governs the main process dynamics, whose impact is observed in the output after a time delay.
The following assumptions are further assumed to the system (\ref{eqn:system_eq_split_incerto}):
\begin{hip}
	\item \label{hip:A-1} The state vector $x_1(t)$ is assumed to be measured in the current time.
	\item \label{hip:A0} The actuator fault $d(t)$ is uniformly bounded such that:
	\begin{align}
		\norm{d(t)}\leq\alpha \label{eqn:NormaFalha}.
	\end{align}
	\item \label{hip:A5} The rate of variation of the delay $\dot{\tau}$ is norm bounded by some known constant $0<\bar{r}<1$
	such that:
	\begin{align}
	\norm{\dot{\tau}}\leq \bar{r}.
	\end{align}
	\item \label{hip:A2} The pair $(A,B)$ is controlable.
	\item \label{hip:A3} The matrix $B$ has full-rank.
		\item \label{hip:A4} The effects of the parametric uncertainties $\delta(x,t)$ are norm bounded by
		\begin{align}
		\norm{\delta(x,t)}\leq\bar{\delta}\norm{x},
		\end{align}
\end{hip}	
where $\bar{\delta}>0$ is a known constant.

This work aims to develop a control strategy capable of stabilizing the system despite the presence of time-varying delays, actuator faults \( d(t) \), and parametric uncertainties. In addition, the approach seeks to reconstruct the fault signal to support system diagnosis.

\section{Predictor and Nonlinear Observer}
In this section, we design an observer-predictor framework aimed at estimating the current value of the unmeasured state component \( x_2 \), as well as reconstructing the disturbance input \( d(t) \) (associated with actuator faults), based on the available output measurement \( y(t) \).

\subsection{Predictor for Measurement Delays} \label{ss:Predx2VOC}
To begin our analysis, let us first consider system (\ref{eqn:system_eq_split_incerto}) without parametric uncertainties ($\delta(x,t)=0$). Under assumption \ref{hip:A-1}, as demonstrated in \cite{Hardy-SBAI2017,Hardy-IJRNC_bayer_mengao}, the $x_2$-dynamics can be characterized as a linear system with known input $x_1$ and delayed output $y=x_2(t-\tau(t))$. Consequently, an estimate $\hat{x}_2$ can be directly obtained from the exact solution of the uncertainty-free system through the variation of constants formula \cite{Artstein-1982,Krstic-2010,Lechappe-2018}:
\smfs\begin{align}
	\hat{x}_2=e^{A_{22}\tau(t)}x_{\tau,2}+\int_{t-\tau(t)}^{t} e^{A{22}(t-\theta)}A_{21}x_1(\theta)d\theta.
	\label{eqn:xh2_VOC}
\end{align}\normalsize

The key benefit of this approach is that, when the system is completely known, the estimate $\hat{x}_2$ serves as an accurate prediction of $x_2$, obtained in finite time from its delayed value $x{\tau,2}$.

In the case of parametric uncertainties, the precise solution for the $x_2$-dynamics would be represented by
\smfs\begin{align}
{x}_2&=e^{A_{22}\tau(t)}x_{\tau,2} +\int_{t-\tau(t)}^{t} e^{A{22}(t-\theta)}\left[A_{21}x_1(\theta)+D_2\delta(x,\theta)\right]d\theta.
\label{eqn:x2_VOC_incerto}
\end{align}\normalsize

In this regard, when utilizing the predictor (\ref{eqn:xh2_VOC}), the prediction would be influenced by a residual error $\tilde{x}_2=x_2-\hat{x}_2$. Conversely, an examination of (\ref{eqn:xh2_VOC}) and (\ref{eqn:x2_VOC_incerto}), we can confirm that
\smfs\begin{align}
\tilde{x}_2&=\int_{t-\tau(t)}^{t} e^{A{22}(t-\theta)}D_2\delta(x,\theta)d\theta.
\label{eqn:x2_VOC_incerto_erro}
\end{align}\normalsize

From the expression (\ref{eqn:x2_VOC_incerto_erro}), it is evident that the duration of the delay significantly impacts the magnitude of this error, which can be upper-bounded by
\smfs\begin{align}
\norm{\tilde{x}_2}\leq \tau(t)\bar{\delta}\norm{e^{A_{22}\tau(t)}D_2} \sup_{\theta \in [t-\tau(t),t]}(\norm{x(\theta)}).
\end{align}\normalsize

\subsection{Nonlinear Super-Twisting Observer}
Using the estimate $\hat{x}_2$, a second-order SMO can be applied with the \textit{Super-Twisting Algorithm} (STA) \cite{Nagesh-2014,Chalanga-2016} to estimate the actuator fault $d(t)$. This is done by defining the estimation error:
\begin{align} \label{meu_pau_bayer_mengo}
	\tilde{x}_1=x_1-\hat{x}_1.
\end{align}
From (\ref{meu_pau_bayer_mengo}), we can design the following nonlinear observer:
\begin{align}
\left\lbrace
	\begin{array}{llll}
		\dot{\hat{x}}_1&=&A_{11}\hat{x}_1+A_{12}\hat{x}_2+B_1 u+\nu\\
		\nu&=&-k_1\frac{\tilde{x}_1}{\norm{\tilde{x}_1}^{1/2}}+k_2\tilde{x}_1+\hat{\xi}\\
		\dot{\hat{\xi}}&=&-k_3\frac{\tilde{x}_1}{\norm{\tilde{x}_1}}-k_4\tilde{x}_1\,.
	\end{array}
\right.
\label{eqn:obsSTA}
\end{align}
Here, $\hat{x}_1$ represents the estimate of the state vector $x_1$, while $\hat{x}_2$ denotes the predicted value of $x_2$, derived from the observer-predictor framework introduced in Section~\ref{ss:Predx2VOC}. It is important to note that the prediction $\hat{x}_2$ is exact when there are no uncertainties in the system parameters or delay. However, in the presence of such uncertainties, the residual error $\tilde{x}_2$ may be nonzero. Consequently, the observer relying on the STA is capable of accommodating these uncertainties as part of the fault signal.

Based on assumptions \ref{hip:A0} and \ref{hip:A4}, the procedure outlined in \cite[eq. (25)]{Nagesh-2014} can be directly utilized to compute the gains $k_1,\cdots,k_4$ and ensure that the sliding mode condition $\dot{\tilde{x}}_1=\tilde{x}_1\equiv 0$ is reached in finite time. Alternatively, similar approaches to those described in \cite{Levant-1993} and \cite{Moreno-2008} can also be applied to determine these gains.

The dynamics of the estimation error $\tilde{x}_1$ can be expressed by incorporating the influence of uncertainties as follows:
\smfs\begin{align}
\left\lbrace
\begin{array}{llll}
\dot{\tilde{x}}_1&=&A_{11}\tilde{x}_1+A_{12}\tilde{x}_2+B_1 d+D_1\delta(x,t)-\nu\\
\nu&=&-k_1\frac{\tilde{x}_1}{\norm{\tilde{x}_1}^{1/2}}+k_2\tilde{x}_1+\hat{\xi}\\
\dot{\hat{\xi}}&=&-k_3\frac{\tilde{x}_1}{\norm{\tilde{x}_1}}-k_4\tilde{x}_1\,.
\end{array}
\right.
\label{eqn:DinErrSTA}
\end{align}\normalsize

After the sliding mode is reached, we can write
\smf\begin{align}
\left\lbrace
\begin{array}{llll}
0&=&0+A_{12}\tilde{x}_2+B_1 d+D_1 \delta(x,t)-\nu\\
\nu&=&\hat{\xi}, \quad
\dot{\hat{\xi}}=0\,.
\end{array}
\right.
\label{eqn:DinErrSTA_deslizamento}
\end{align}\normalsize

The first equation in (\ref{eqn:DinErrSTA_deslizamento}) can be rearranged to obtain:
\smf\begin{align}
\hat{\xi}&=B_1 d+D_1\delta(x,t)+A_{12}\tilde{x}_2,
\label{eqn:definehatxi}
\end{align}\normalsize
which serves as an estimate of the combined effects of faults and disturbances on $x_1$. It is worth noting that, in the absence of uncertainties, the signal $\hat{\xi}$ can be directly employed for reconstructing the fault $d(t)$.

\subsection{Observer-Predictor Combination}
At this point, the predictor (\ref{eqn:xh2_VOC}) and the observer (\ref{eqn:obsSTA}) can be integrated into a unified observer-predictor framework designed for uncertain systems featuring output delays and partially unmeasured states.
\smfs\begin{align}
\left\lbrace
\begin{array}{lll}
		\hat{x}_2&=&e^{A_{22}\tau(t)}x_{\tau,2}+\int_{t-\tau(t)}^{t} e^{A{22}(t-\theta)}A_{21}x_1(\theta)d\theta\\
		\dot{\hat{x}}_1&=&A_{11}\hat{x}_1+A_{12}\hat{x}_2+B_1 u+\nu\\
		\nu&=&-k_1\frac{\tilde{x}_1}{\norm{\tilde{x}_1}^{1/2}}+k_2\tilde{x}_1+\hat{\xi}\\
		\dot{\hat{\xi}}&=&-k_3\frac{\tilde{x}_1}{\norm{\tilde{x}_1}}-k_4\tilde{x}_1\,,
\end{array}
\right.
\label{eqn:observer_full_version}
\end{align}\normalsize
where $\tilde{x}_1 =x_1-\hat{x}_1$. This integrated observer-predictor structure enables finite-time estimation of the state $x_2$ and the disturbance term combined with parametric uncertainties, denoted by $\hat{\xi}$. When there are no uncertainties, the actuator fault $d(t)$ can also be directly reconstructed from the signal $\hat{\xi}$, in accordance with assumption \ref{hip:A3}.

\section{Sliding Mode Control}
After computing the estimates $\hat{x}_2$ and $\hat{\xi}$, and based on assumption \ref{hip:A2}, a SMC can be formulated to ensure asymptotic stabilization while maintaining robustness against model uncertainties and external disturbances.

\subsection{Control Law}
First of all, let us define the following state vector
\smfs\begin{align}
	\bar{x}=\vechx,
\end{align}\normalsize
which includes the measured state $x_1$ and the predicted value $\hat{x}_2$, provided by the predictor (\ref{eqn:xh2_VOC}), both evaluated at the current time instant.

Next, we define the sliding variable
\smfs\begin{align}
s(t)=S\bar{x}=\left[\begin{array}{cc}I&S_2\end{array}\right]\vechx,
\end{align}\normalsize
where $S$ is an appropriate matrix to be designed such that $SB$ is not singular.

Equation (\ref{eqn:xh2_VOC}) yields the expression for $\dot{\hat{x}}_2$
\smfs\begin{align}
\dot{\hat{x}}_2&=(1-\dot{\tau})e^{A_{22}\tau}\dot{x}_{\tau,2}+\dot{\tau}A_{22}e^{A_{22}}{x}_{\tau,2} +A_{21}x_1-(1-\dot{\tau})e^{A_{22}\tau}A_{21}{x}_{\tau,1}
+A_{22}\int_{t-\tau}^{t} e^{A_{22}(t-\theta)}A_{21}{x}_1(\theta)d\theta.
\end{align}\normalsize

By substituting equations (\ref{eqn:system_eq_split_incerto}) and (\ref{eqn:xh2_VOC}) accordingly, it is possible to verify that
\smfs\begin{align}
\dot{\hat{x}}_2&= A_{21}x_1 + A_{22}\hat{x}_2 + (1-\dot{\tau})e^{A_{22}\tau}D_2 \delta(x,t-\tau).
\label{eqn:dotxhat2}
\end{align}\normalsize
To simplify the notation, we define
\smfs\begin{align}
\zeta_1(x,t)&:=D_1 \delta(x,t),\\
\zeta_2(x,t)&:=(1-\dot{\tau})e^{A_{22}\tau}D_2 \delta(x,t-\tau).
\label{eqn:def_zeta}
\end{align}\normalsize

Let us examine the linear transformation
\smf\begin{align}
T:=\left[\begin{array}{cc}I&S_2\\0&I\end{array}\right], 
\end{align}\normalsize
such that  
\smfs\begin{align}
\left[\begin{array}{c}s\\\hat{x}_2\end{array}\right]=T\vechx\,.
\end{align}\normalsize

Therefore, the dynamic system can be expressed in terms of the new coordinates as follows:
\smfs\begin{align}
\left\lbrace
\begin{array}{lll}
\left[\begin{array}{c}\dot{s}\\\dot{\hat{x}}_2\end{array}\right]&=&\underbrace{\splitAbar}_{\bar{A}} \left[\begin{array}{c}s\\\hat{x}_2\end{array}\right]+ \splitB (u+d) 
+\left[\begin{array}{c}\zeta_1(x,t)+S_2\zeta_2(x,t)\\\zeta_2(x,t)\end{array}\right]\,,\\
y(t)&=& x_{\tau,2}\,,
\end{array}
\right.
\label{eqn:dinamica_SMC}
\end{align}\normalsize
where $\bar{A}_{11}=A_{11}+S_2A_{21}$, $\bar{A}_{12}=S_2A_{22}+A_{12}-\bar{A}_{11}S_2$,
$\bar{A}_{21}=A_{21}$ and $\bar{A}_{22}=A_{22}-A_{21}S_2$.

To design the matrix $S$, we first consider equation (\ref{eqn:dinamica_SMC}) for $\dot{\hat{x}}_2$, assuming no parametric uncertainties, i.e., when $\zeta_2(x,t)=0$. During the sliding mode, $s=0$ and
\smfs\begin{align}
\dot{\hat{x}}_2&=  \bar{A}_{22}\hat{x}_2\,.
\end{align}\normalsize
Therefore, during the sliding mode, the dynamics of $\hat{x}_2$ are determined by the eigenvalues of $\bar{A}{22} = A_{22} - A_{21} S_2$. Moreover, assumption \ref{hip:A2} guarantees that the pair $(A_{22}, A_{21})$ is controllable (see \cite[Prop. 3.3]{Spurgeon-book}), meaning that any pole-placement method can be used to select $S_2$ so that $\bar{A}_{22}$ is Hurwitz.

Consider the control law
\smfs\begin{align}
u=u_d+u_{\rm nom}+u_{\rm SM},
\end{align}\normalsize
where $u_d$ is the term designed for disturbance rejection. In this case, $u_{\rm nom}$ represents the nominal component of the control signal, while $u_{\rm SM}$ is the switching term, which is used to drive the sliding variable $s$ to the sliding mode where $s \equiv 0$.

It is important to note that the $x_1$-dynamics in (\ref{eqn:system_eq_split_incerto}) are
\smfs\begin{align}
\dot{x}_1=A_{11}x_1+A_{12}x_2 + B_1(u+d) + \zeta_1(x,t).
\label{eqn:x_1}
\end{align}\normalsize
From (\ref{eqn:definehatxi}), we can define:
\smfs\begin{align}
u_d=-B_1^\dagger \hat{\xi}=-d-B_1^\dagger D_1\delta(x,t)-B_1^\dagger A_{12}\tilde{x}_2,
\label{eqn:ud}
\end{align}\normalsize
where $B_1^\dagger$ is the Moore-Penrose pseudo-inverse of $B_1$, as expressed by $B_1^\dagger= (B_1^T B_1)^{-1} B_1^T$ \cite[page 122]{Strang-book}.

Upon insertion of (\ref{eqn:ud}) into (\ref{eqn:x_1}), the following relationship is verified:
\smfs\begin{align}
\dot{x}_1=A_{11}x_1+A_{12}\hat{x}_2 + B_1[u_{\rm nom}+u_{\rm SM}].
\label{eqn:x_1+ud}
\end{align}\normalsize
Given that $\tilde{x}_2 = x_2 - \hat{x}_2$, it is evident that the control input $u_d$ is capable of compensating for the disturbance $d$, as well as the matched components of the parametric uncertainties, specifically $A_{12} \tilde{x}_2$ and $\zeta_1(x,t)$.

Based on the structure of $u_d$, the nominal control input $u_{\rm nom}$ can be derived using the equivalent control method~\cite{Shtessel-book}, by solving the condition $\dot{s} = 0$ for the nominal (uncertainty-free) system, while taking into account that the dynamics of $x_1$ are already expressed as in equation~(\ref{eqn:x_1+ud}). Assuming $u_{\rm SM} = 0$, it follows that we can write:
\smfs\begin{align}
0=SA\bar{x}+SB(u_{\rm nom}). 
\end{align}\normalsize
%
%
Given that 
$SB$ is nonsingular by design, it follows that
\smfs\begin{align}
u_{\rm nom}=-(SB)^{-1}SA\bar{x}. 
\end{align}\normalsize

To conclude the control design, the switching component is introduced to enforce the sliding condition, ensuring that the system trajectories reach and remain on the sliding manifold.
\smfs\begin{align}
u_{\rm SM}=-\rho (SB)^{-1} {\rm sign}(s),
\end{align}\normalsize
where ${\rm sign}(s) = \left[{\rm sign}(s_1)  \;\cdots \;{\rm sign}(s_{n-p}) \right]^T$ and the parameter $\rho > 0$ must be selected large enough to dominate the residual uncertainties and disturbances that are not mitigated by the control components $u_d$ and $u_{\rm nom}$.

Consequently, the complete form of the control law is expressed as:
\smfs\begin{align}
u=-B_1^\dagger \hat{\xi}-(SB)^{-1}SA\bar{x}-\rho (SB)^{-1} {\rm sign}(s).
\label{eqn:leidecontrole}
\end{align}\normalsize

\subsection{Stability Analysis}
By connecting (\ref{eqn:leidecontrole}) into (\ref{eqn:dinamica_SMC}), we get
\smfs\begin{align}
\left\lbrace
\begin{array}{lll}
\dot{s}&=&S_2\zeta_2(x,t)-\rho {\rm sign}(s)\\
\dot{\hat{x}}_2&=&\bar{A}_{21}s+\bar{A}_{22}\hat{x}_2+\zeta_2(x,t)\\
y&=&{x}_{\tau,2}=\hat{x}_{\tau,2}+\tilde{x}_{\tau,2}.\\
\end{array}
\right.
\label{eqn:malha_fechada}
\end{align}\normalsize

Based on the closed-loop system described in (\ref{eqn:malha_fechada}), one can confirm that the sliding mode condition $s \equiv 0$ is attained in finite time, as stated in the following theorem.

\medskip
\begin{theo}
	\label{teo:mododeslizante}
	The control law (\ref{eqn:leidecontrole}) applied to (\ref{eqn:system_eq_split_incerto}), under the assumptions
	\ref{hip:A-1} to \ref{hip:A4}, with 
	\smfs\begin{align}
	\rho\geq\phi\bar{\delta}(1+\bar{r}) \norm{e^{A_{22}\tau}D_2} \norm{x}+\eta\,,
	\label{eqn:expressao_rho}
	\end{align}\normalsize
	$\eta>0$ and $\phi>1$ being arbitrary constants, guarantees the existence of the sliding mode $s=S\bar{x}\equiv0$ in finite time.
\end{theo}
\medskip
\proof 
Consider the Lyapunov function $V=s^Ts/2$. By computing the time derivative of \( V \) along the trajectories of the system in~(\ref{eqn:malha_fechada}), and employing the definition of \( \zeta_2(x,t) \) from~(\ref{eqn:def_zeta}), we get
\smfs\begin{align}
\dot{V}&=s^TS_2\zeta_2(x,t)-\rho s^T{\rm sign}(s),\\
\dot{V}&=(1-\dot{\tau})s^T e^{A_{22}\tau}D_2 \delta(x,t-\tau)-\rho s^T{\rm sign}(s).\label{eqn:Vdot}
\end{align}\normalsize
Notice that 
\smfs\begin{align}
 s^T{\rm sign}(s)&=\norm{s}_1\geq\norm{s}_2, \\
 -\rho s^T{\rm sign}(s)&\leq-\rho\norm{s}_2.
\end{align}\normalsize
Since \( \dot{V} \) includes a delayed term, we will apply Razumikhin's Theorem~\cite[Theorem 1.4]{Gu-book} in the subsequent analysis.  
To begin with, recall that the condition \( \|x(t+\theta)\| < \phi \|x(t)\| \) holds for some \( \phi > 1 \) and \( \theta \in [-\tau(t), 0] \). Hence,
\smfs\begin{align}
\norm{\delta(x,t+\theta)}&\leq\bar{\delta}\norm{x(t+\theta)}\leq \phi\bar{\delta}\norm{x(t)},\\
\norm{\delta(x,t-\tau)}&\leq\phi\bar{\delta}\norm{x(t)}. \label{eqn:condRzmk}
\end{align}\normalsize
Therefore, employing \ref{hip:A5}, we have
\smfs\begin{align}
\dot{V}&\leq -\norm{s} \left[\rho -\phi\bar{\delta}(1+\bar{r}) \norm{e^{A_{22}\tau}D_2} \norm{x(t)}\right].
\label{eqn:Vdot_final}
\end{align}\normalsize
If
\smfs\begin{align}
\rho\geq \phi\bar{\delta}(1+\bar{r}) \norm{e^{A_{22}\tau}D_2} \norm{x}+\eta,
\end{align}\normalsize
then
\smfs\begin{align}
\dot{V}&\leq -\eta \norm{s}<0,
\end{align}\normalsize
demonstrating that $s$ is globally asymptotically stable. Furthermore, since
\smfs\begin{align}
\dot{V}&\leq -\eta \sqrt{V},
\end{align}\normalsize
%
one concludes the sliding surface condition $s\equiv0$ is achieved within finite time. $\hfill \square$
\endproof
\medskip


The expression in (\ref{eqn:Vdot}) reveals a term dependent on the delay variation rate $\dot{\tau}$, originating from the perturbation signal $\zeta_2(x,t)$. The derivation of the upper bound in (\ref{eqn:Vdot_final}) requires $\dot{\tau}$
  to be bounded - a condition guaranteed by Assumption \ref{hip:A5}. Notably, in the absence of parametric uncertainties (when $\delta(x,t)=0$ in (\ref{eqn:Vdot})), this assumption could be relaxed as the $\dot{V}$ expression would no longer contain delay-dependent terms.

After the sliding mode is established, we can confirm that the closed-loop system (\ref{eqn:malha_fechada}) is asymptotically stable through the subsequent corollary.

\medskip
\begin{crl} 
	Given the same assumptions \ref{hip:A-1} to \ref{hip:A4} from Theorem \ref{teo:mododeslizante}, once the sliding mode is reached, the closed-loop system (\ref{eqn:malha_fechada}) is asymptotically stable if there exists a matrix $P_2>0$ such that
 	$P_2\bar{A}_{22}+\bar{A}_{22}^TP_2=-I$ and
	\smfs\begin{align}
	\bar{\delta}< \frac{1}{2\phi\lambda_{\max}(P_2)(1+\bar{r}) \left(\sqrt{1+\norm{S_2}^2}\right) \norm{e^{A_{22}\tau}D_2}}, 
	\label{eqn:Bound_incerteza}
	\end{align}\normalsize
where \( \phi > 1 \) is a chosen constant and \( \lambda_{\max}(P_2) \) denotes the largest eigenvalue of the matrix \( P_2 \). \label{crl:Corolario}
\end{crl}
\medskip
\proof
Once the sliding mode is established, we have $s\equiv0$ and the resulting system dynamics, obtained from (\ref{eqn:malha_fechada}), 
is given by:
\smfs\begin{align}
\dot{\hat{x}}_2=\bar{A}_{22}\hat{x}_2+\zeta_2(x,t).
\label{eqn:din_modo_deslizante}
\end{align}\normalsize

Reminding that $\bar{A}_{22}$ is Hurwitz by design, therefore, We are able to determine a symmetric positive definite matrix $P_2$ such that $P_2\bar{A}_{22}+\bar{A}_{22}^TP_2=-I$. Subsequently, we define the Lyapunov function as follows: $V_2(\hat{x}_2)=\hat{x}_2^TP_2\hat{x}_2$. The time derivative along the trajectories of (\ref{eqn:din_modo_deslizante}) is expressed as:
\smfs\begin{align}
\dot{V}_2(\hat{x}_2)&=-\hat{x}_2^T\hat{x}_2+2\hat{x}_2^TP_2\zeta_2(x,t)\\
&\leq -\norm{\hat{x}_2}^2+2\lambda_{\max}(P_2)\norm{\hat{x}_2}\norm{\zeta_2(x,t)}\nonumber\\
&= -\lambda_{\max}(P_2)\norm{\hat{x}_2}\left[\mu\norm{\hat{x}_2}-2\norm{\zeta_2(x,t)}\right],\nonumber
\label{eqn:V2dot}
\end{align}\normalsize
where $\mu:=1/\lambda_{\max}(P_2)$.

From (\ref{eqn:def_zeta}), we know that
\smfs\begin{align}
\zeta_2(x,t)&=(1-\dot{\tau})e^{A_{22}\tau}D_2 \delta(x,t-\tau).
\end{align}\normalsize
Since it depends on a delayed variable, we will apply Razumikhin's Theorem once more.
In the sliding mode, it holds that $x_1=-S_2\hat{x}_2$, such that  $\norm{x}\leq\left(\sqrt{1+\norm{S_2}^2}\right)\norm{\hat{x}_2}$. Thus, from (\ref{eqn:condRzmk}), we get
\smfs\begin{align}
\norm{\delta(x,t-\tau)}&\leq\phi\bar{\delta}\left(\sqrt{1+\norm{S_2}^2}\right)\norm{\hat{x}_2}.
\label{eqn:bound_Delta}
\end{align}\normalsize

In this regard, by substituting (\ref{eqn:bound_Delta}) into (\ref{eqn:V2dot}), we obtain
\smfs\begin{align}
\dot{V}_2(\hat{x}_2)  \leq &-\lambda_{\max}(P_2)\norm{\hat{x}_2}^2\left[\mu - \beta_1\right],\\
\beta_1:=&2\phi\bar{\delta}(1+\bar{r}) \left(\sqrt{1+\norm{S_2}^2}\right) \norm{e^{A_{22}\tau}D_2}. 
\end{align}\normalsize

Whenever $\beta_1 < \mu$, it follows that $\dot{V}_2(\hat{x}_2) < 0$. According to Razumikhin's Theorem, this guarantees the asymptotic stability of the $\hat{x}_2$-dynamics under the condition that  
\smfs\begin{align}
2\phi\bar{\delta}(1+\bar{r}) \left(\sqrt{1+\norm{S_2}^2}\right) \norm{e^{A_{22}\tau}D_2}< \frac{1}{\lambda_{\max}(P_2)},
\end{align}\normalsize
resulting in the condition (\ref{eqn:Bound_incerteza}). $\hfill \square$
\endproof
\medskip

Corollary \ref{crl:Corolario} establishes the convergence $\hat{x}_2 \to 0$ as $t \to \infty$. During sliding mode operation, the relation $x_1 = -S_2\hat{x}_2$ implies that $x_1 \to 0$ as well. Furthermore, the inequality $\norm{x} \leq \left(\sqrt{1+\norm{S_2}^2}\right)\norm{\hat{x}_2}$ holds, and assumption \ref{hip:A4} guarantees that perturbation effects asymptotically vanish. Consequently, equation (\ref{eqn:x2_VOC_incerto_erro}) yields $\tilde{x}_2 \to 0$, which in turn implies $x_2 \to 0$. This proves that the closed-loop system achieves stabilization despite actuator faults and parametric uncertainties, while satisfying the bound condition (\ref{eqn:Bound_incerteza}).

While assumption \ref{hip:A4} enables exact stabilization at the origin, this requirement can be relaxed by considering a more general upper bound for the parametric uncertainties, such as $\norm{\delta(x,t)} \leq \bar{\delta}_1\norm{x} + \bar{\delta}_2$. Under this weaker condition, the system state would converge to an ultimately bounded region around the origin, rather than achieving asymptotic convergence to the origin itself.

\section {Numerical Example}
To demonstrate the effectiveness of the proposed approach, we first examine the nominal system (\ref{eqn:system_eq})--(\ref{eqn:system_eq_split_incerto}) without parametric uncertainties ($\delta(x,y) \equiv 0$), where
\smf\begin{align}
A=
\left[
	\begin{array}{cc}
	 -1&1\\-3&1
	\end{array}
\right], \quad
B=
\left[
\begin{array}{c}
1\\0
\end{array}
\right], \quad
C=
\left[
\begin{array}{cc}
1&0
\end{array}
\right].
\end{align}\normalsize
The actuator fault was modeled as $d(t) = 0.1\sin(2t) + 0.2\cos(3t)$, with a time-varying output delay given by $\tau(t) = 0.4 + 0.1\sin(t)$ s. Observer gains were selected following the methodology in \cite[eq. (25)]{Nagesh-2014}, resulting in $k_1 = k_3 = 5$ and $k_2 = k_4 = 2$. For the control implementation, we employed the control law (\ref{eqn:leidecontrole}) with parameters $S = [1\; -5]$ and a constant gain $\rho = 2$ for simplicity.

\begin{figure} [htb]
\centering
\includegraphics[scale=0.55]{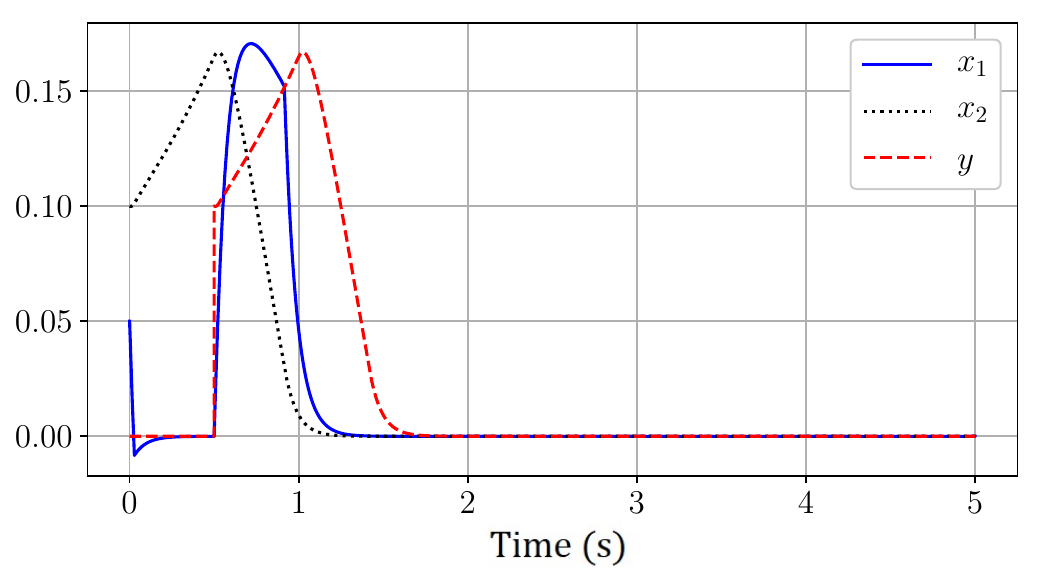}
\caption{\underline{Without Parametric Uncertainties}: state variables $x_1(t)$ (in blue), $x_2(t)$ (in black) and the output signal $y(t)$ (in red).}
\label{fig:teste61}
\includegraphics[scale=0.55]{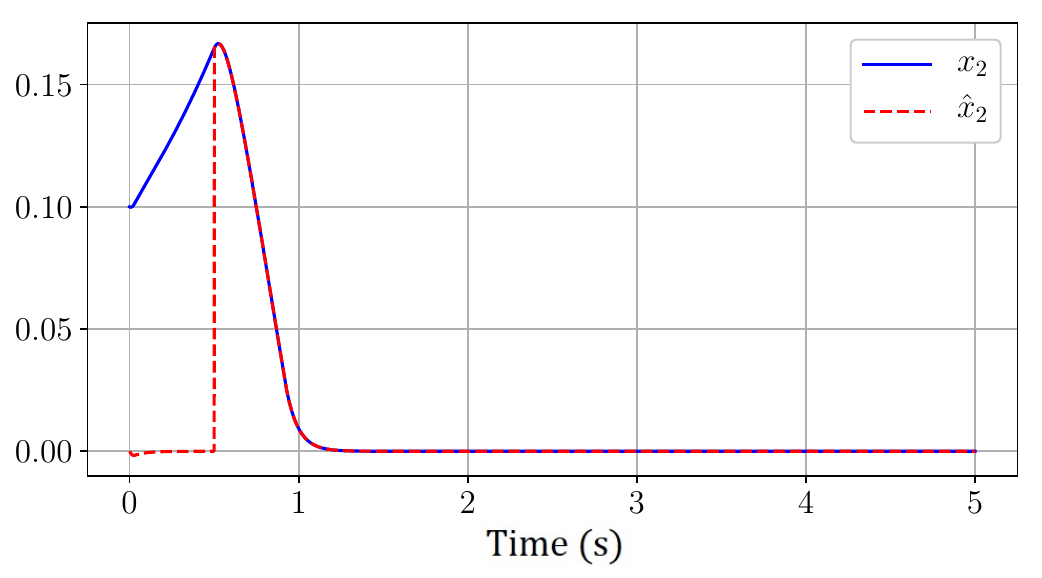} [htb]
\caption{\underline{Without Parametric Uncertainties}: state variable $x_2(t)$ (in blue) and its estimate $\hat{x}_2(t)$ (in red), obtained with the finite-time predictor (\ref{eqn:xh2_VOC}).}
\label{fig:teste63}
\end{figure}

The system response is presented in Figure \ref{fig:teste61}, displaying the state variables $x_1$ (blue) and $x_2$ (black), along with the output signal $y$ (red). The proposed controller successfully stabilizes the plant using both the measured state $x_1$ and the observer-predictor estimates $\hat{\xi}$, $\hat{x}_2$ generated by (\ref{eqn:observer_full_version}). 

\begin{figure} [htb]
\centering
\includegraphics[scale=0.55]{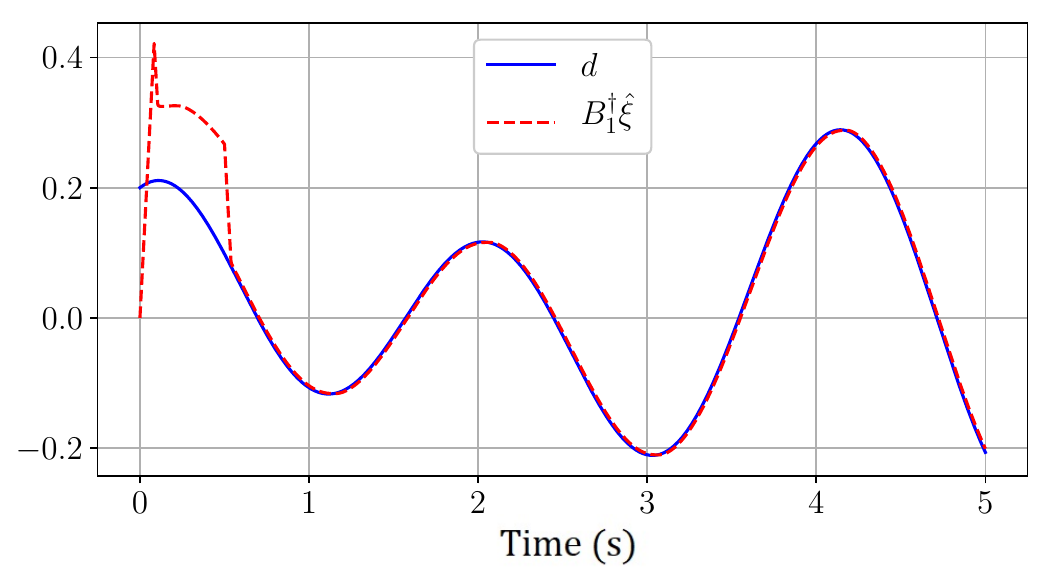}
\caption{\underline{Without Parametric Uncertainties}: disturbance $d(t)$ (in blue) and its estimate $\hat{d}(t)\!=\!B_1^\dagger\hat{\xi}(t)$ (in red), obtained with the finite-time observer (\ref{eqn:obsSTA}).}
\label{fig:teste64}
\includegraphics[scale=0.55]{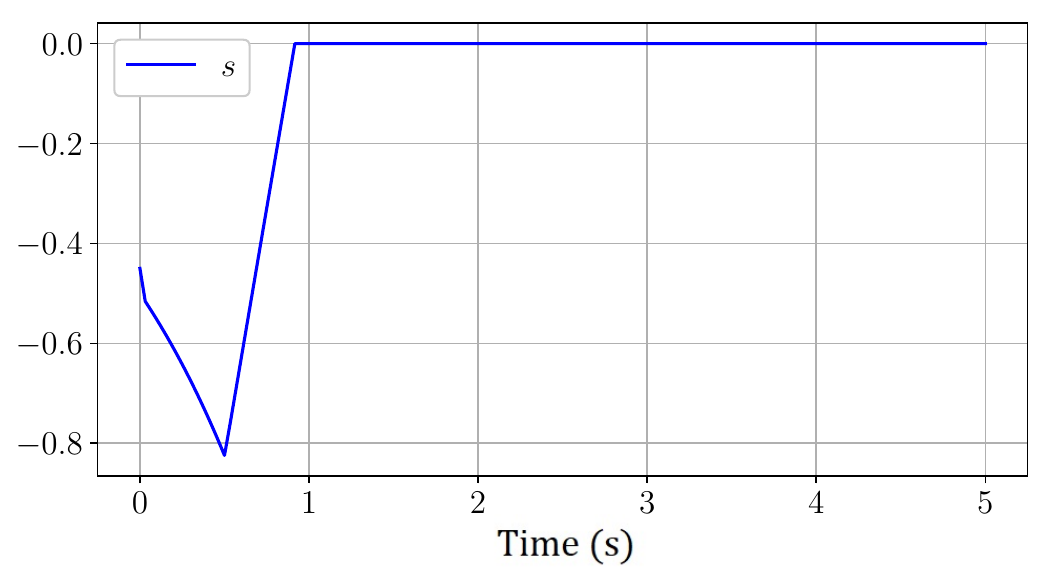}
\caption{\underline{Without Parametric Uncertainties}: sliding variable $s(t)$ reaching zero in finite time.}
\label{fig:teste65}
\includegraphics[scale=0.55]{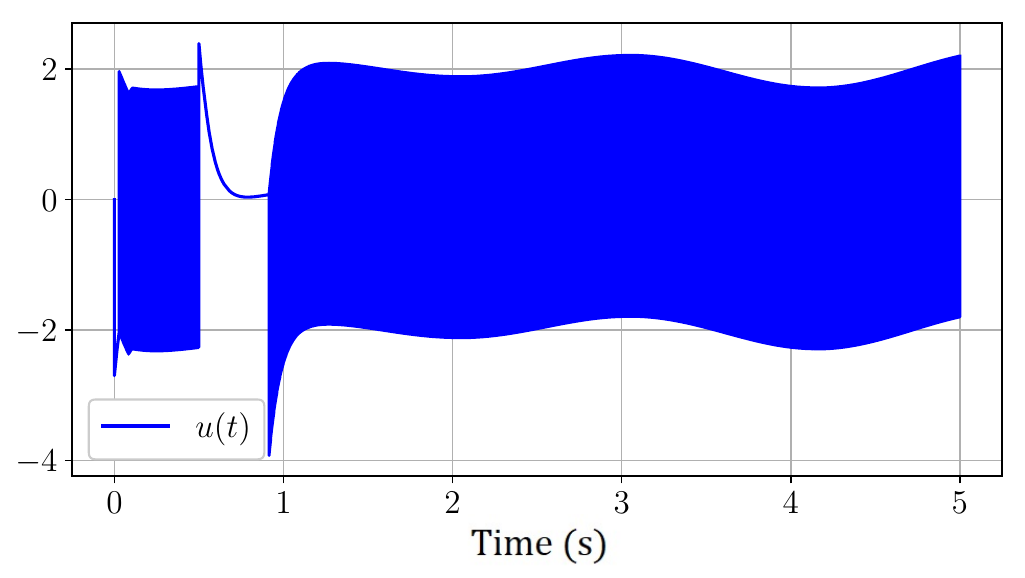}
\caption{\underline{Without Parametric Uncertainties}: control signal $u(t)$ highlighting the occurrence of the sliding mode.}
\label{fig:teste66}
\end{figure}

Figure \ref{fig:teste63} demonstrates the finite-time convergence of the predictor (\ref{eqn:xh2_VOC}), comparing the actual state $x_2$ (blue) with its estimate $\hat{x}_2$ (red). Similarly, Figure \ref{fig:teste64} presents the actuator fault $d$ alongside its finite-time estimate $\hat{d} = B_1^\dagger \hat{\xi}$, where $B_1^\dagger = (B_1^T B_1)^{-1} B_1^T$ denotes the left pseudo-inverse of $B_1$ \cite[page 122]{Strang-book}, computed by the observer (\ref{eqn:obsSTA}).

The sliding surface behavior and control effort are shown in Figures \ref{fig:teste65} and \ref{fig:teste66}, respectively, illustrating the system's performance during stabilization.

\begin{figure} [htb]
\centering
\includegraphics[scale=0.5]{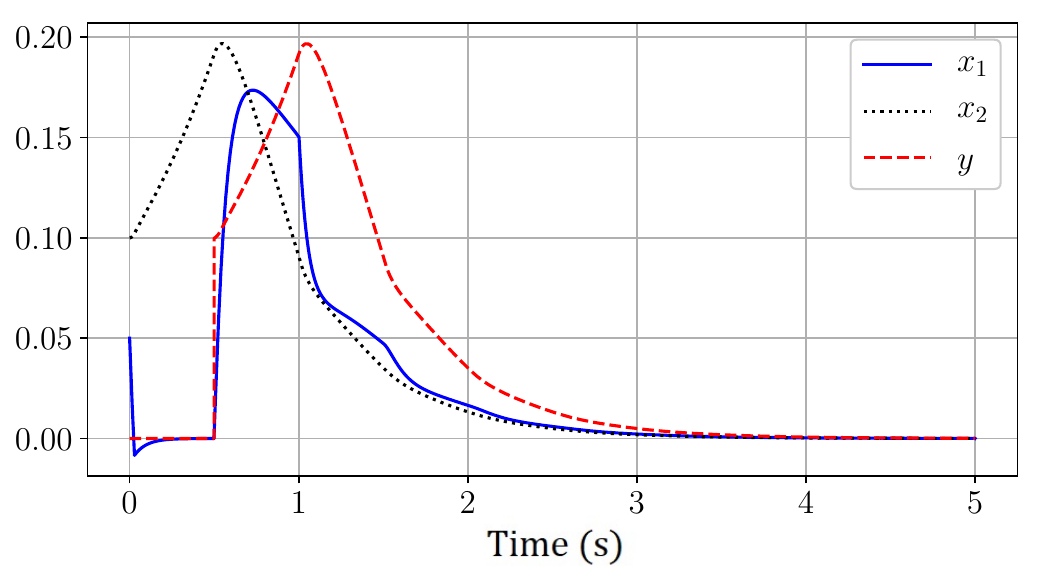}
\caption{\underline{With Parametric Uncertainties}: state variables $x_1(t)$ (in blue), $x_2(t)$ (in black) and the output signal $y(t)$ (in red).}
\label{fig:teste61i}
\includegraphics[scale=0.55]{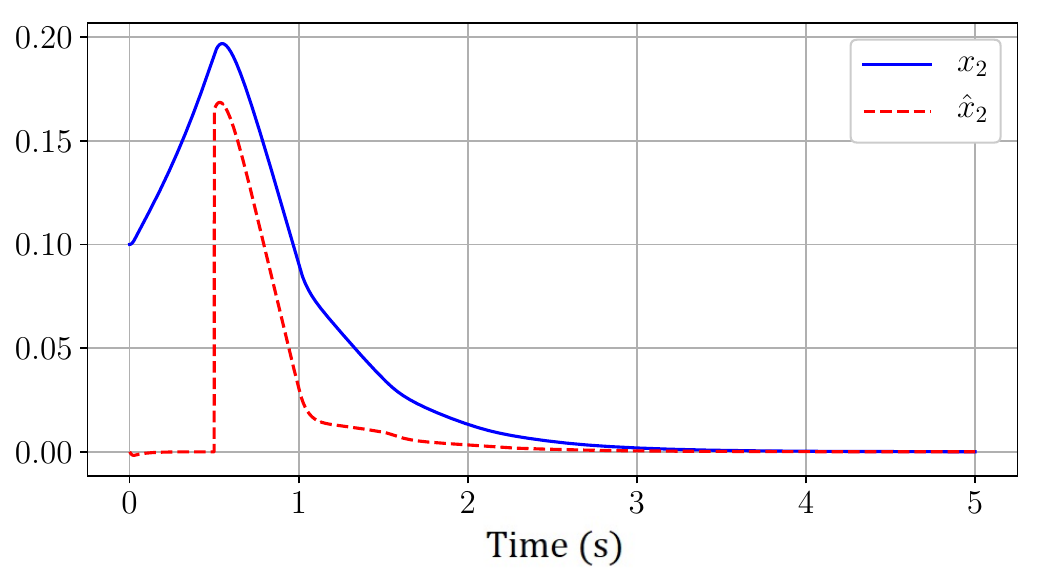}
\caption{\underline{With Parametric Uncertainties}: state variable $x_2(t)$ (in blue) and its estimate $\hat{x}_2(t)$ (in red), obtained with the finite-time predictor (\ref{eqn:xh2_VOC}).}
\label{fig:teste63i}
\includegraphics[scale=0.55]{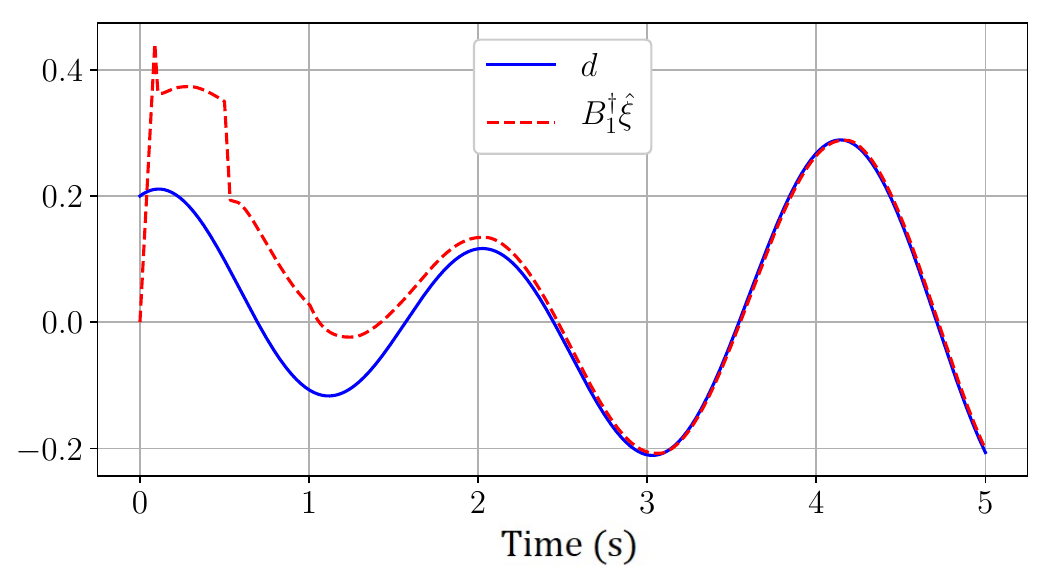}
\caption{\underline{With Parametric Uncertainties}: disturbance $d(t)$ (in blue) and its estimate $\hat{d}(t)\!=\!B_1^\dagger\hat{\xi}(t)$
(in red), obtained by means of (\ref{eqn:definehatxi}). The estimate tries to accommodate the uncertainties effects in the prediction error $\tilde{x}_2(t)$.}
\label{fig:teste64i}
\end{figure}

To assess the robustness of the proposed method under parametric uncertainties, we selected
\smfs\begin{align}
D=\left[
\begin{array}{cc}
0.4&0.4\\0.4&0.4
\end{array}
\right], \quad \delta(x,t)=x.
\end{align}\normalsize
For this second case, the control gain was redefined as $\rho=5$ to handle the parametric uncertainties and maintain the stability of the closed-loop system—without this adjustment, the system would become unstable (curves omitted).
As seen in Figure \ref{fig:teste61i}, the system is still stabilized even in the presence of such uncertainties.
Figure \ref{fig:teste63i} shows how these uncertainties impact the estimation accuracy of $\hat{x}_2$.
Meanwhile, Figure \ref{fig:teste64i} illustrates that the nonlinear Super-Twisting observer is capable of reconstructing the disturbance $d$, thus compensating for the negative effects introduced by the uncertainties.


\newpage
\section{Conclusions}
This paper proposed a novel strategy for stabilizing linear time-invariant systems affected by time-varying measurement delays, parametric uncertainties, and nonlinear matched disturbances—such as those representing actuator faults. An additional complexity considered was the partial unavailability of the system state vector.
The proposed solution combines an open-loop predictor with a nonlinear observer grounded in the Super-Twisting Algorithm.
Differently from previous studies, this class of observers enabled the estimation of a broader category of disturbance signals while preserving robustness against parametric uncertainties and delay effects.
A key benefit of this structure is that, upon reaching sliding mode in finite time via the predictor-observer approach, the closed-loop system could still be stabilized, at least asymptotically.

Future investigation lies in the expansion of proposed design and analysis for different control problems with unknown control direction as a particular class of faults and pursuing single- or multi-agent optimization via extremum or Nash equilibrium seeking rather than only stabilization, as
considered in references \cite{Oliveira_controldirection,TRO_book2022,Oliveira_jogos}. Adaptive sliding mode control \cite{Oliveira-Cunha-Hsu-Springer-2017} and estimation \cite{Rodrigues-Oliveira-IJC-2018} are also fruitful scenarios for further developments.

\end{document}